# Teaching mathematics with a different philosophy

# Part 2: Calculus without Limits

C. K. Raju
School of Mathematical Sciences
Universiti Sains Malaysia
11700, Penang, Malaysia
ckr@ckraju.net

**Abstract:** The example of the calculus is used to explain how simple, practical math was made enormously complex by imposing on it the Western religiously-colored notion of mathematics as "perfect". We describe a pedagogical experiment to make math easy by teaching "calculus without limits" using the new realistic philosophy of zeroism, different from Platonic idealism or formalist metaphysics. Despite its demonstrated advantages, it is being resisted because of the existing colonial hangover.

## *1. Introduction*

Part 1 explained how Western mathematics originated in mathesis and religious beliefs about the soul. Hence, mathematics was first banned by the church, and later reinterpreted, in a theologically-correct way, as a "universal" metaphysics. This post-Crusade reinterpretation (based on the myth of "Euclid") was not historically valid, for the *Elements* did use empirical proofs. Eventually, Hilbert and Russell eliminated empirical proofs (in the *Elements* and mathematics) and made mathematics fully metaphysical. Even so, this metaphysics is *not* universal, but has a variety of biases, as was pointed out. It has nil practical utility. In contrast, most math of practical value originated in the non-West with a different epistemology,[1] which permitted empirical proofs (which do not diminish practical value in any way). While the West adopted this math for its practical value, it tried to force-fit it into its religious beliefs about math: first that math must be "perfect" (since it incorporates eternal truths), and, second, that this perfection could only be achieved through metaphysics, since the empirical world was considered imperfect. (Hence, the idea that math *must* be metaphysical.)

Present-day learning difficulties in math reflect the historical difficulties that arose in this way, because the West imposed a religiously biased Western metaphysics on practical non-Western mathematics.

An example might make matters clearer.

## *2. The case of the calculus*

Contrary to the false history that the calculus originated with Newton and Leibniz, it has now been firmly established that the calculus originated in India.[2] The process started in the 5th c., with Aryabhata's attempts to calculate precise trigonometric values. That precision was needed for *practical* reasons for accurate astronomical models, and the calendar—essential for monsoon-dependent agriculture—as also for navigation.[3] (Agriculture and overseas trade were the two key sources of wealth in India then.)

In the 16th c., Europe was poor, and European dreams of wealth rested on overseas trade. This presupposed good navigation technique. However, Europeans had peculiar difficulties with navigation.

First their reliance on charts and straight lines led to the problem of loxodromes. (Indo-Arabic navigators did *not* have this problem, for they used celestial navigation,[4] and the 7th c. Bhaskara had already mentioned[5] the objection that the sphericity of the earth must be taken into account while calculating latitude and longitude.) Second, Europeans could not use techniques such as those of Bhaskara, to determine longitude, because Columbus grossly underestimated the size of the earth, leading to the 1501 Portuguese ban on carrying globes aboard ships. They could not even determine latitude in daytime since their religious calendar (Julian calendar) was off by 11 days in the 16th c. The (European) navigation problem was recognized as the key scientific challenge in Europe then. It remained so for centuries, and many governments offered huge rewards for its solution. (The last was the British prize, for a method to determine longitude at sea; this was legislated in 1712, and partly given away around 1760.)

Accurate trigonometric values were needed to calculate loxodromes (to make the Mercator chart), and also to determine latitude and longitude. The most precise trigonometric values then available anywhere in the world were in Indian timekeeping texts (precise to 9 decimal places[6]), which texts[7] were found in the vicinity of Cochin. The first Roman Catholic mission also happened to be established in Cochin, in 1500, and by the mid-16th c., the Jesuits took over a well-established Christian college in Cochin which they used as a base to acquire Indian texts, translate them, and despatch them back to Europe, in Toledo mode. The authors of some of the key Indian astronomy texts, such as Sankara Variyar who wrote the *Yuktidipika*, even shared a common patron with the Portuguese in the Raja of Cochin. Clavius published those precise Indian trigonometric values under his name in 1607,[8] though he did not know enough trigonometry to calculate the size of the earth! To summarise the new history, *the calculus originated in India for its practical value* (for agriculture and navigation) *and was brought to Europe also for its practical value* (for navigation, and especially the problems of determining the three ells—loxodromes, latitude, and longitude[9]).

However, those precise trigonometric values were calculated by Indian mathematicians using infinite series expansions (today called "Taylor's" expansion, "Leibniz" series, etc.), and sophisticated techniques[10] to sum infinite series. These techniques were not comprehended by European mathematicians (who were, then, still struggling at the level of decimal fractions introduced by Stevin, only in 1582). The key difficulty was with the notion of infinite sums, as in the non-terminating, non-recurring decimal expansion for the number $\pi$. The notion of infinity brought religious beliefs prominently into play.

Thus, the number $\pi$ could also be finitely understood as the ratio of the circumference to the diameter of a circle. However, Descartes objected that the ratio of the length of a curved line with a straight line was beyond the human mind.[11] Now here is such a simple thing—a child can use a string to measure the length of a curved line, and straighten it to compare it to that of a straight line. That was how mathematics was taught in India[12] since the days of the *sulba sutra*—or "aphorisms on the string"—but a major Western philosopher asserted this to be impossible! Descartes' difficulty arose from his religious beliefs about mathematics as "perfect". Hence, he naively imagined that "rigorously" obtaining the circumference of a circle required one to break up the circumference into straight-line segments, and physically sum up the lengths of the segments, leading to the infinite ("Leibniz") series. Descartes thought such an infinite sum could only be done by God. Stopping the sum at a finite stage would mean neglecting a small quantity; though irrelevant for *all* practical applications, it would make mathematics "imperfect", hence not mathematics at all. Galileo concurred, and *hence* left matters to his student Cavalieri, to avoid the risk of disrepute.

Newton thought he had resolved these difficulties of Descartes and Galileo, about infinity, by a clever appeal to God! He needed the notion of time derivative for his second "law" of motion. He thought this notion of derivative could be made "rigorous" by his doctrine of fluxions. This required that time itself must "flow" ("smoothly", or "equably").[13] Now, while things may flow *in* time, the slightest thought shows that this idea that time itself flows is meaningless, and has long been recognized as such.[14] Nevertheless, Newton thought mathematics was the "perfect" language in which God had written the "laws" of nature. He admitted that time could not be properly measured by physical phenomena which were "imperfect". But he postulated a perfect, "absolute, true, and mathematical time", which "flows equably" but "without relation to anything external".[15] Each adjective, "absolute", "true", "mathematical", shows that Newton thought time was metaphysical and known only to God, and if anyone still had a doubt, he added the last clause "without relation to anything external". People often quote Newton on this without understanding that making time metaphysical was the weakest point of his physics, and his physics failed exactly for that reason, and had to be replaced by relativity.[16] This shows how metaphysical considerations regarding mathematics have impeded science.

These Western difficulties with the Indian calculus continued with Berkeley's objections[17] to the illogical procedures about infinitesimals/fluxions used by Leibniz and Newton, to which Newton's supporters had no serious answer.[18] Eventually, Dedekind brought in formal reals, $R$, but this required the metaphysical manipulation of infinity enabled by Cantor's set theory. That, in turn, was suspect and was formalised only in the 1930's. Formal set theory is so difficult that only a few mathematicians bother to learn it—the head of the math department of an IIT could not even state the formal definition of a set, when publicly challenged to do so by this author. Naturally, students who learn the "new math", which begins with set theory, find math difficult. Ironically, this formalisation led to the belated realization that calculus can also be done over "non-Archimedean" fields, larger than $R$, such as the field of rational functions used by Indian mathematicians who treated rational functions much like ordinary fractions. In such a field, limits are not unique, unless one discards infinitesimals, and that procedure (exactly what Indian mathematicians adopted) is equivalent to limits by order counting.

Note that all this theologising about the "perfect" way to handle infinity, which went through curved lines, fluxions, formal reals, limits, and sets, and has returned to infinitesimals, has added not an iota to the *practical* value of calculus. Aryabhata's numerical method[19] (equivalent to what is today called "Euler's method" for ordinary differential equations) is still adequate for *all* practical applications of calculus to Newtonian physics. (Of course, the method can be and has been improved.) Naturally, engineering students still ask today "what is the the point of doing limits?", and the teachers have no answer except to recite the magical word "rigor". As we have already seen, this claim of rigor has no substance, but incorporates merely an unreasonable demand that mathematics must conform to a particular, religiously-biased metaphysics.

This theological Western view of math was globalised by the *political* force of colonialism. It was stabilised by Macaulay's well known intervention with the education system, and the continued support for it is readily understood on Huntington's doctrine of soft-power.[20] And this way of teaching math continues to be uncritically followed to this day even after independence. This is the first attempt to try to re-examine and critically re-evaluate the Western philosophy of math[21] and suggest an alternative to European ethnomathematics.

The new philosophy proposed by this author has now[22] been renamed "zeroism",[23] to emphasize that it is being used for its practical value, and does not depend upon (the interpretation of) any Buddhist texts about *sunyavada*. A key idea is that of mathematics as an adjunct physical theory. Another key idea is that, like infinitesimals, small numbers may be neglected, as in a computer calculation, but on the new grounds that *ideal representations are erroneous*, for they can *never* be achieved in reality (which is continuously changing). (Exactly what constitutes a discardable "small" number, or a "practical infinitesimal", is decided by the context,[24] as with formal infinitesimals or order-counting.) This is the antithesis of the Western view that mathematics being "ideal" must be "perfect", and that only metaphysical postulates for manipulating infinity (as in set theory), laid down by authoritative Western mathematicians, are reliable, and all else is erroneous. As this debate between realism and idealism is an old one, we will not go into further details here.

## 3. The experiment

If the learning difficulties with math arise from the theological complexities that the West has woven into math, and if that math was universalised by the political force of colonialism, then the natural remedy is to decolonise math by dispensing with those theological complexities, which anyway add nothing to the practical value of math. And the fact is that the mass of students today learn math for its practical value. These considerations led to the new course on calculus without limits, which aims to teach calculus using (a) zeroism, (b) computers, and (c) by following the actual historical trajectory of the development of calculus as concerning the numerical solution of ordinary differential equations.

A recent experiment, over the last couple of years, has tested the feasibility and desirability of teaching this new course. The experiment involved five groups till now, one at the Central University of Tibetan Studies, Sarnath,[25] and four groups at the School of Mathematical Sciences, Universiti Sains Malaysia (USM). The group sizes varied from 6 to 35. The 4 groups at the USM consisted of one group of post-graduate math students, one group of undergraduate pure math students, one group of undergraduate applied math students, and one group of non-math students. The availability of four groups at USM allowed one to test *separately* the various claims about the course as follows.

Among the advantages claimed for the new course are the following.
(1) The new philosophy makes the calculus easier to understand. Thus, any calculus student today can parrot off that $\frac{d}{dx}e^x = e^x$, but, few (even among IIT students) have even a rough idea of the *definition* of $e^x$. In the new approach, functions are rigorously defined as the solution of differential equations. Thus $y = e^x$ is the rigorously defined and explained as the solution of $y' = y$ with the condition that $y(0) = 1$. Students are quickly able to calculate the values of the function, plot it, and analyse it in various ways using software such as this author's CALCODE.
(2) The new approach allows students to apply calculus to advanced problems. For example, the motion of the simple pendulum involves the Jacobian elliptic functions; consequently, most people learn only the simplified theory of the simple pendulum, and often confound it with simple harmonic motion.[26] With the new approach, the Jacobian elliptic function sn *x* is not particularly more complicated than sin *x*, so students can study the variation in the time period of the simple pendulum with its amplitude. Similar remarks apply to the problem of ballistics with air resistance, or why a heavier cricket ball can be thrown further than a tennis ball.
(3) The new approach makes calculus so easy that even non-math students can master it in a short while.

With all five groups, there was a pre-test and a post-test. The post-tests included 6 to 10 problems drawn at random (using a pseudo-random-number generator) from a published question bank of 12[th] standard calculus. Those questions involved symbolic calculation of symbolically complex derivatives and integrals, which is what students are expected to master in current calculus courses. This was included in the post-test only to demonstrate that teaching those skills is completely pointless today, when it can be done in a jiffy using open-source symbolic manipulation programs such as MAXIMA, That is *not* the same as reliance on a calculator for arithmetic sums: in daily life one occasionally needs to do arithmetic sums in one's head, but one never needs to do in one's head any complicated integrals or derivatives (involving only elementary functions!). The other aspects of the post-test included solving and analysing the solution of ordinary differential equations, since that is the at the heart of applications of the calculus. Also included were some non-elementary elliptic integrals.

The performance on the post-test was uniformly good even for the non-math students (with at least middle-school level math). Thus, claims (2) and (3) were validated. However, the pre-test revealed that even the post-graduate students were not well-versed with the philosophy of formal mathematics, and were not comfortable with advanced mathematical notions such as the Schwartz derivative. As such it was not possible to test whether they found the deeper aspects of the new philosophy easier than the existing philosophy of formalism.

The whole approach can be extended to several variables and partial differential equations in an obvious way. But that is a future agenda.

## *4. The dimension of hegemony*

Identifying the difficulties with math learning, and proposing a solution, does represent a major advance. But there are difficulties in implementing the solution. Various stakeholders (such as students afraid of math, or their parents) are *never* consulted to decide what sort of math to teach. Even scientists and engineers are rarely consulted regarding what sort of math ought to be taught to them. However, if all decisions regarding the math curriculum are left solely to math "experts", there is an obvious conflict of interests: for these experts were brought up on the older tradition of formal mathematics, and rejecting formal mathematics may well make their past work valueless.

Public discussion is one way to ensure that the interests of millions of students are not disadvantaged, and that scientific and educational activities relate to public interest.[27] Such discussions would be particularly welcome given the other sensitive issue in the present case: namely that imposing a religiously biased metaphysics on millions of students is not only unethical, it is unconstitutional under the Indian constitution which guarantees secularism, or under any other constitution which does not permit a Western religious bias.